\documentclass[11pt]{amsart}

\usepackage[T1]{fontenc}
\usepackage[utf8]{inputenc}
\usepackage{lmodern}
\usepackage{amsmath,amssymb,amsthm,mathtools}
\usepackage{enumitem}
\usepackage[a4paper,margin=2.75cm]{geometry}
\usepackage{microtype}
\usepackage[hidelinks]{hyperref}
\usepackage{aliascnt}
\usepackage{color}
\usepackage{xcolor}

\usepackage[nameinlink,capitalize,noabbrev]{cleveref}
\hypersetup{
  pdftitle={Transversally Cantor-like limit sets for meromorphic connections},
}

\setlist[enumerate]{label=(\roman*),leftmargin=2.3em}
\setlist[itemize]{leftmargin=2em}

\newtheorem{theorem}{Theorem}[section]

\newaliascnt{proposition}{theorem}
\newtheorem{proposition}[proposition]{Proposition}
\aliascntresetthe{proposition}

\newaliascnt{lemma}{theorem}
\newtheorem{lemma}[lemma]{Lemma}
\aliascntresetthe{lemma}

\newaliascnt{corollary}{theorem}
\newtheorem{corollary}[corollary]{Corollary}
\aliascntresetthe{corollary}

\newaliascnt{claim}{theorem}

\aliascntresetthe{claim}

\newaliascnt{conjecture}{theorem}

\aliascntresetthe{conjecture}

\theoremstyle{definition}
\newaliascnt{definition}{theorem}
\newtheorem{definition}[definition]{Definition}
\aliascntresetthe{definition}

\newaliascnt{example}{theorem}
\newtheorem{example}[example]{Example}
\aliascntresetthe{example}

\newaliascnt{question}{theorem}

\aliascntresetthe{question}

\theoremstyle{remark}
\newaliascnt{remark}{theorem}
\newtheorem{remark}[remark]{Remark}
\aliascntresetthe{remark}

\crefname{theorem}{theorem}{theorems}
\Crefname{theorem}{Theorem}{Theorems}
\crefname{proposition}{proposition}{propositions}
\Crefname{proposition}{Proposition}{Propositions}
\crefname{lemma}{lemma}{lemmas}
\Crefname{lemma}{Lemma}{Lemmas}
\crefname{corollary}{corollary}{corollaries}
\Crefname{corollary}{Corollary}{Corollaries}
\crefname{definition}{definition}{definitions}
\Crefname{definition}{Definition}{Definitions}
\crefname{remark}{remark}{remarks}
\Crefname{remark}{Remark}{Remarks}
\crefname{question}{question}{questions}
\Crefname{question}{Question}{Questions}
\crefname{conjecture}{conjecture}{conjectures}
\Crefname{conjecture}{Conjecture}{Conjectures}

\newcommand{\C}{\mathbb C}
\newcommand{\R}{\mathbb R}

\newcommand{\Sone}{\mathbb S^1}
\newcommand{\supp}{\operatorname{supp}}
\newcommand{\Sing}{\operatorname{Sing}}
\newcommand{\Res}{\operatorname{Res}}
\newcommand{\Hol}{\operatorname{Hol}}

\newcommand{\Int}{\operatorname{Int}}

\title[Cantor-like limit sets]{Transversally Cantor-like limit sets for meromorphic connections}

\author{Kholikova Madina}
\address{National University of Uzbekistan named after Mirzo Ulugbek, 4 University Street,
100174 Tashkent, Uzbekistan}
\email{khalikovamadina9@gmail.com}

\author{Eralieva Ozoda}
\address{National University of Uzbekistan named after Mirzo Ulugbek, 4 University Street,
100174 Tashkent, Uzbekistan}
\email{ozodaeralieva23023@gmail.com}


\begin{document}

\begin{abstract}
We study transversally Cantor-like $\omega$-limit sets of eventually simple geodesics for meromorphic connections on compact Riemann surfaces. We show that their minimal invariant subsets are precisely the poles and establish the sharp bound $\operatorname{Re}\operatorname{Res}_p\nabla\geq -1/2$ for Fuchsian poles in nontrivial limit sets. Moreover, we rule out Cantor-like recurrence for connections with real periods, regardless of their rotational holonomy, thereby refining the Poincar\'e--Bendixson classification.
\end{abstract}

\maketitle
\enlargethispage{2pt}

\section{Introduction}
 Let $S$ be a compact connected Riemann surface. A \emph{meromorphic
connection} $\nabla$ on its tangent bundle $TS$ is a $\C$-linear operator
which assigns to each local meromorphic vector field $X$ a meromorphic
$1$-form with values in $TS$ and satisfies the Leibniz rule
$$
 \nabla(fX)=df\otimes X+f\nabla X
$$
for every local meromorphic function $f$. In a local coordinate $z$ it is
represented by a meromorphic $1$-form $\eta$, determined by
$\nabla(\partial_z)=\eta\otimes\partial_z$. The poles of these local
representations form a finite set $\Sing(\nabla)$; we write
$S^0=S\setminus\Sing(\nabla)$. A pole $p$ is \emph{Fuchsian} if $\eta$ has
a simple pole there, and the connection is Fuchsian if all its poles are
Fuchsian. Its residue $\Res_p\nabla:=\Res_p\eta$ is independent of the
coordinate.

A \emph{geodesic} is a nonconstant smooth real curve in $S^0$ satisfying
$\nabla_{\sigma'}\sigma'=0$. If $\eta=f(z)\,dz$, this equation becomes
$z''+f(z)(z')^2=0$. For a geodesic $\sigma:[0,\varepsilon)\to S^0$ which is maximal
in the forward direction, where $0<\varepsilon\leq+\infty$, its $\omega$-limit set is
$$
\omega(\sigma)=\bigcap_{0<t<\varepsilon}\overline{\sigma([t,\varepsilon))},
$$
where closures are taken in $S$. Thus poles can belong to the limit set
even though the geodesic is contained in $S^0$. We call $\sigma$
\emph{eventually simple} if its restriction to $[t_0,\varepsilon)$ is injective for
some $t_0<\varepsilon$. Deleting an initial segment does not change its
$\omega$-limit set.

Meromorphic connections provide a common setting for several classes of
geometric dynamical systems. In suitable coordinates their geodesics are
Euclidean straight lines, connecting complex differential equations with
affine geometry and surface topology. Abate--Tovena studied their limiting
behaviour on the Riemann sphere and its relation to homogeneous
holomorphic vector fields \cite{AT}; Abate--Bianchi established
Poincar\'e--Bendixson results on arbitrary compact Riemann surfaces
\cite{AB}. Connections with meromorphic differentials, singular flat
metrics and affine surfaces are developed in \cite{NST,R}.

 One recurrent possibility requires particular care: a simple geodesic
may be contained in its own $\omega$-limit set, while that set has empty
interior and Cantor intersections with transverse curves. This
possibility is included in the classification of Abate--Rakhimov
\cite[Appendix~A]{AR}. Its closure in the compact surface can contain
poles, so the regular transverse dynamics and the invariant subsets of
the whole limit set must both be considered.

 We recall the meanings of these terms before stating our results. The
associated-line-field construction of \cite{AB} assigns to every regular
point of the limit set of a simple geodesic a distinguished tangent line,
obtained by taking limits of its tangent lines. The corresponding local
integral curves contained in the limit set are geodesics. A closed subset is \emph{invariant} if it
contains the entire maximal distinguished geodesic through each of its regular points;
poles are regarded as singular invariant points. A nonempty closed
invariant set is \emph{minimal} if it contains no proper nonempty closed
invariant subset. A \emph{closed regular geodesic} is a geodesic whose
support is a smooth closed curve contained in $S^0$.

 We call $W=\omega(\sigma)$ \emph{transversally Cantor-like} if, after
deleting an initial segment, $\sigma$ is injective and contained in $W$,
$\Int W=\varnothing$, and $I\cap W$ is a Cantor set for every sufficiently
short geodesic transversal $I$ through the remaining geodesic. Here a
transversal is a regular geodesic segment meeting the distinguished
tangent lines transversely; it is taken compact with endpoints outside
$W$. A Cantor set means a nonempty compact perfect totally disconnected
set. These conditions imply $W=\overline{\sigma([t_0,\varepsilon))}$; see
\Cref{def:cantor}.

Our first result concerns arbitrary meromorphic connections and
describes the singular part of such a limit set.

\begin{theorem}\label{thm:core-intro}
Let $\nabla$ be a meromorphic connection on a compact Riemann surface $S$, and let $W$ be a transversally Cantor-like $\omega$-limit set of an eventually simple geodesic. Put
$$
 P_W:=W\cap\Sing(\nabla).
$$
Then:
\begin{enumerate}
\item $P_W$ is finite and nonempty;
\item the minimal closed invariant subsets of $W$ are exactly the singletons $\{p\}$ with $p\in P_W$. Every Fuchsian pole $p\in P_W$ satisfies $\operatorname{Re}\Res_p\nabla\geq-\frac{1}{2}$;
\item $W$ contains no closed regular geodesic.

\end{enumerate}
\end{theorem}

 The description of the singular core follows from the minimal-set
classification of Abate--Bianchi, based on Hounie's theorem
\cite{AB,Hou}. The residue bound also uses the local geometry of
geodesics. A Fuchsian pole with $\operatorname{Re}\Res_p\nabla\leq-1$
can belong only to a singleton $\omega$-limit set \cite{R}. In the
remaining range $-1<\operatorname{Re}\rho<-1/2$, where
$\rho=\Res_p\nabla$, a local coordinate gives
$\eta=\rho\,dz/z$. Writing $a=1+\rho$, the multivalued developing map
$J(z)=z^a/a$ sends geodesics to straight lines. We prove an exact
self-intersection criterion for their inverse images: when
$0<\operatorname{Re}a<1/2$, a sufficiently deep passage entering and
leaving a fixed small neighborhood of the pole intersects itself.
This contradicts eventual simplicity, including when $\rho$ is not real.

 Our second result concerns connections with real periods. On $S^0$ one
can choose local coordinates $w_\alpha$ in which the connection form
vanishes. Their transition maps have the form
$$
 w_\beta=a_{\alpha\beta}w_\alpha+c_{\alpha\beta},
 \qquad a_{\alpha\beta}\in\C^*.
$$
Continuation around loops defines the \emph{linear holonomy}
$\Hol_\nabla:\pi_1(S^0)\to\C^*$. The connection has \emph{real periods}
if its linear holonomy is contained in
$\Sone=\{u\in\C:|u|=1\}$. In this case the affine coordinates can be
normalized so that $|a_{\alpha\beta}|=1$. Their Euclidean metrics agree,
giving a flat metric on $S^0$; in the Fuchsian case this is an adapted
singular flat metric on $S$, whose geodesics are those of $\nabla$
\cite{R}. The rotational holonomy need not be finite.

\begin{theorem}\label{thm:flat-intro}
Let $\nabla$ be a Fuchsian meromorphic connection with real periods on a compact Riemann surface. Then no transversally Cantor-like $\omega$-limit set can occur.
\end{theorem}

 The proof uses a short transversal contained entirely in the regular
part of the surface. Consecutive crossings cut the simple geodesic into
disjoint return arcs. Keeping track of the crossing direction and the
relative topological type gives finitely many return families; arcs of
the same type bound nonsingular Euclidean quadrilaterals. The resulting
monotone return correspondences transport complementary gaps of the
recurrent crossing set, with only finitely many exceptions. A cyclic
gap chain would produce a closed regular geodesic, excluded by
\Cref{thm:core-intro}, so an infinite chain of distinct gaps must exist.
Their lengths tend to zero. Unitary parallel transport preserves the
difference between the directions at the two endpoints, and continuity
therefore forces this difference to vanish. The boundary geodesics are
then parallel, and their positive perpendicular separation is preserved
along the chain, contradicting the shrinking gap lengths.

 Combining \Cref{thm:flat-intro} with \cite{AR} gives the strengthened
Poincar\'e--Bendixson classification in \Cref{cor:PB}. The real-period
condition is global: real residues alone do not imply unitary holonomy.
The genus-two example discussed at the end of the paper illustrates
this distinction using the Disco surface of \cite{BFG}.

 The paper is organized as follows. In \cref{sec:prelim} we recall the
preliminary notions and the results on associated line fields.
\Cref{sec:core} contains the local self-intersection argument and the
proof of \Cref{thm:core-intro}. In \cref{sec:real period} we construct
the regular return bands and prove \Cref{thm:flat-intro}.
\Cref{sec:PB} gives the strengthened Poincar\'e--Bendixson theorem and
discusses the role of nonunitary holonomy.

\subsection*{Acknowledgement.}

The authors would like to express their sincere gratitude to their supervisor, Professor Karim Rakhimov, for valuable guidance, insightful suggestions, and continuous support throughout this research. The authors also gratefully acknowledge the support provided by Grant No FL-9524115114, which made this research possible.

\subsection*{AI Declaration.}

The main idea of this paper emerged from a year-long study of the literature, focusing on the non-isometric nature of affine interval exchange maps and the existence of wandering intervals. ChatGPT 6 was used assist in searching and verifying mathematical arguments, and improving the clarity and accuracy of the manuscript. The authors take full responsibility for the content of the paper.

\section{Preliminary notions}\label{sec:prelim}

In this section we repeat some definitions and theorems from \cite{AB}, \cite{AR} and \cite{R}.
\subsection{Meromorphic connections and geodesics}

\label{sec:connections} \begin{definition}\label{def:meromorphic connection} 
 Let $TS$ be the tangent bundle of a Riemann surface $S$. A \emph{holomorphic connection} on $TS$ is a $\C-$linear map  $\nabla:\mathcal TS\longrightarrow
  \Omega_S^1\otimes\mathcal TS$ satisfying the Leibniz rule
$$\nabla(fs)=df\otimes s+f\nabla s$$
for all $f \in \mathcal{O}_S$ and $s \in \mathcal TS$, where $\mathcal TS$ denotes the sheaf of germs of holomorphic sections of $TS$, while $\mathcal{O}_S$ is the sheaf of germs holomorphic functions on $S$ and $ \Omega_S^1$ is the sheaf of germs of holomorphic 1-forms on $S$.

Let $(U_\alpha , z_\alpha)$ be a local chart for $S$, and $\nabla$ a holomorphic connection on $S$. It is easy to see that there exists a holomorphic 1-form $\eta_\alpha \in \Omega_S^1(U_\alpha)$ on $U_\alpha$ such that
$$ \nabla(\partial_\alpha)=\eta_\alpha\otimes\partial_\alpha,$$ 
where $\partial_\alpha:= \frac{\partial}{\partial z_\alpha}$ is the induced local generator of $TS$ over $U_\alpha$.
\end{definition}

\begin{definition}
    We say that $\eta_\alpha$ is the \emph{local representation} of $\nabla$ on $U_\alpha$.
    \vspace{\baselineskip}
    
    Now, for given local representations $\{ \eta_\alpha\}$ we look for a condition which guarantees the existence of holomorphic connection $\nabla$. Let $\{\xi_{\alpha \beta} \}$ be the cocycle representing the cohomology class $\xi \in H^1(S, \mathcal{O}^*)$ of $TS$ (hence $\xi_{\alpha \beta}=\frac{\partial z_\alpha}{\partial z_\beta}$); over $U_\alpha \cap U_\beta$ we have 
    $$\partial_\beta=\partial_\alpha \xi_{\alpha \beta}$$
    and thus 
    $$\nabla (\partial_\beta) = \nabla(\partial_\alpha \xi_{\alpha \beta}) \Leftrightarrow \eta_\beta \otimes\partial_\beta = \xi_{\alpha \beta} \eta_\alpha \otimes\partial_\alpha + d\xi_{\alpha \beta} \otimes \partial_\alpha$$ 
    if and only if 
    $$\eta_\beta=\eta_\alpha + \frac{d\xi_{\alpha\beta}}{\xi_{\alpha\beta}}.$$
\end{definition}

\begin{definition}
    Let $\nabla$ be a holomorphic connection on a Riemann surface $S$. We say that a Hermitian metric $g$ on $TS$ is \emph{adapted} to $\nabla$ if $\nabla g\equiv 0$, that is if 
    $$X(g(R,T))=g(\nabla_X R,T) +g(R, \nabla_{\overline{X}}T)$$
    and 
    $$\overline{X}(g(R,T))=g(\nabla_{\overline{X}}R,T)+g(R,\nabla_X T)$$
    for all sections $R,T$ of $TS$, and every vector field $X$ on $S$.
\end{definition}

\begin{definition} 
A \emph{meromorphic connection} on the tangent bundle $TS$ of a Riemann surface $S$ is a $\C-$linear map  $\nabla:\mathcal M_{ TS}\longrightarrow
 \mathcal M_S^1\otimes\mathcal M_{TS}$ satisfying the Leibniz rule
$$\nabla(fs)=df\otimes s+f\nabla s$$
for all $f \in \mathcal M_S$ and $s \in \mathcal M_{TS}$, where ${ \mathcal M_{TS}}$ denotes the sheaf of germs of meromorphic sections of $TS$, while $\mathcal M_S$ is the sheaf of germs meromorphic functions on $S$ and $\mathcal M_S^1$ is the sheaf of germs of meromorphic 1-forms on $S$.
    
\end{definition}

Let $(U_\alpha , z_\alpha)$ be a local chart for $S$, and $\nabla$ a meromorphic connection on $S$. Then there exists $\eta_\alpha \in \mathcal M_S^1(U_\alpha)$ on $U_\alpha$ such that
$$ \nabla(\partial_\alpha)=\eta_\alpha\otimes\partial_\alpha,$$ 
where $\partial_\alpha:= \frac{\partial}{\partial z_\alpha}$ is the induced local generator of $TS$ over $U_\alpha$. As in the case of holomorphic connection we have 
$$\eta_\beta=\eta_\alpha + \frac{d\xi_{\alpha\beta}}{\xi_{\alpha\beta}} $$ on $U_\alpha \cap U_\beta$, where $\xi_{\alpha \beta}:=\frac{\partial z_\alpha}{\partial z_\beta}$. In particular, if every representation is holomorphic then $\nabla$ is a holomorphic connection. We say $p\in S$ is a $pole$ for meromorphic connection $\nabla$ if $p$ is a pole of $\eta_\alpha$ any local chart $U_\alpha$ at $p$. If $\Sigma$ is the set of poles of $\nabla$, then $\nabla$ is a holomorphic connection on $S^0=S\setminus\Sigma$. So we can define a notion of geodesic for $\nabla$ on $S^0$ as in holomorphic connection.
\begin{definition}
    A \emph{geodesic} for a meromorphic connection $\nabla$ on $TS$ is a nonconstant smooth real curve $\sigma:I\to S^0$, with $I\subseteq \R$ an interval, such that $\nabla_{\sigma'}\sigma'=0.$
\end{definition}

Let $\{ (U_\alpha, z_\alpha) \}$ be an atlas for $S$ and $\sigma:I\to U_\alpha$, with $I\subseteq \R$ an interval, a smooth curve. Then $\sigma$ is a geodesic for meromorphic connection $\nabla$ if and only if 
$$(z_\alpha \circ \sigma)''+(f_\alpha\circ \sigma)(z_\alpha \circ \sigma)'^2\equiv0.$$
where $ \eta_\alpha=f_\alpha\,dz_\alpha$ is the local representation of $\nabla$ on $U_\alpha$.

Every point of $S^0$ has a simply connected neighborhood on which
$\eta=dF$.  If $w$ is a primitive of $e^Fdz$, then the representation of
$\nabla$ in the coordinate $w$ is zero.  We call $(U,w)$ a
$\nabla$-\emph{chart}.  On an overlap of two such charts,
\begin{equation}\label{eq:affine-transition}
 w_2=a_{12}w_1+b_{12},
 \qquad a_{12}\in\C^*,\quad b_{12}\in\C.
\end{equation}

Thus $S^0$ carries a complex affine structure, and geodesics are Euclidean
straight lines with affine parametrization in every $\nabla$-chart.

A pole is \emph{Fuchsian} when it is a simple pole of a local representation. Its residue is denoted by $\Res_p\nabla$. The residue of a meromorphic connection at a pole is independent of the choice of local charts.

On a simply connected regular chart, choose a primitive $K$ of $\eta$ and then a primitive $J$ of $e^K$. The map $J$ is a local developing map, or local isometry, for the connection: it sends geodesics to Euclidean straight segments. Two such maps differ by a complex affine transformation.

The linear parts of the changes between developing coordinates define the monodromy representation
$$
 \Hol_\nabla:\pi_1(S^0)\longrightarrow\C^*.
$$
We say that $\nabla$ has \emph{real periods} when
$$
 \Hol_\nabla\bigl(\pi_1(S^0)\bigr)\subset\Sone.
$$
For a Fuchsian connection this is equivalent to the existence of an adapted singular flat metric \cite{AR,AT,R}.

\subsection{Associated line fields and minimal sets}

We recall the main topological tool from \cite{AB}. A smooth line field with singularities on $S$ is a smooth section of the real projectivized tangent bundle on the complement of a closed singular set. A closed set is invariant if it is a union of maximal leaves and singular points.

\begin{definition}
A smooth line field with singularities on a Riemann surface $S$
is a smooth section of the real projectivized tangent bundle
$ \mathbb{P}_{\R}(TS)$ on the complement of a closed singular set.
\end{definition}

\begin{theorem}[Abate and Bianchi \cite{AB}\label{thm:associated}]
Let $\sigma$ be a maximal simple geodesic of a meromorphic connection, maximal in both time directions. There exists a smooth line field $\Lambda$ with singularities on $S$ such that:
\begin{enumerate}
\item $\sigma$ is an integral curve of $\Lambda$;
\item in a neighborhood of $\supp\sigma$, the singularities of $\Lambda$ are exactly the poles of $\nabla$;
\item $\omega(\sigma)$ is $\Lambda$-invariant;
\item the restriction of $\Lambda$ to $\omega(\sigma)$ is uniquely determined by $\sigma$.
\end{enumerate}
 For one-sided maximal simple geodesics, we use directly the minimal-set classification stated below.
\end{theorem}

Hounie's theorem says that a minimal set of a smooth singular line field on a compact connected surface is a singularity, a closed integral curve, or the whole surface; in the last case the surface is a torus and the line field is topologically equivalent to an irrational linear line field \cite{Hou}. Applied to \Cref{thm:associated}, this gives the following form.

\begin{theorem}[Abate and Bianchi \cite{AB}\label{thm:minimal-AB}]
Let $\sigma:[0,\varepsilon)\to S^0$ be a maximal simple geodesic and let $W=\omega(\sigma)$. Every minimal closed invariant subset of $W$ is one of the following:
\begin{enumerate}
\item a pole of $\nabla$;
\item the support of a closed geodesic, in which case $W$ equals that support;
\item all of $S$, in which case $S$ is a torus and $\nabla$ is holomorphic.
\end{enumerate}
\end{theorem}

\begin{definition}\label{def:cantor}
Let $\sigma:[0,\varepsilon)\to S^0$ be a maximal eventually simple geodesic. Its $\omega$-limit set $W$ is called a \emph{transversally Cantor-like $\omega$-limit set} if, after deleting an initial segment of $\sigma$,
\begin{enumerate}
\item $\sigma$ is simple and $\supp\sigma\subset W$;
\item $\Int W=\varnothing$;
\item for every sufficiently short simple geodesic transversal $I$ through a point of $\supp\sigma$, taken compact with endpoints outside $W$, the set $I\cap W$ is perfect and totally disconnected.
\end{enumerate}
\end{definition}

Since $W$ is closed and is always contained in $\overline{\supp\sigma}$, the condition $\supp\sigma\subset W$ implies
$W=\overline{\supp\sigma}.$
This is precisely the empty-interior subcase discussed in \cite{AR}.

Throughout this section $\nabla$ is an arbitrary meromorphic connection on a compact Riemann surface, $\sigma$ is a simple geodesic satisfying \Cref{def:cantor}, and
$$
 W:=\omega(\sigma)=\overline{\supp\sigma}.
$$
The recurrent simple tail extends to a maximal two-sided simple
distinguished geodesic contained in $W$: a transverse crossing would
contradict the uniqueness of the limiting line, while a tangential
return would make the geodesic closed. We apply
\Cref{thm:associated} to this extension and denote the resulting
line field by $\Lambda$.

\begin{lemma}\label{lem:continuum}
The set $W$ is a compact connected nondegenerate continuum.
\end{lemma}

\begin{proof}
For $t<\varepsilon$, set
$$
 K_t:=\overline{\sigma([t,\varepsilon))}.
$$
Every $K_t$ is compact and connected, and the family is nested as $t$ increases. Moreover,
$$
 W=\bigcap_{t<\varepsilon}K_t.
$$
The intersection of a nested family of nonempty compact connected sets is connected. The set is nondegenerate because a singleton has no Cantor transverse section.
\end{proof}

\section{Proof of Theorem 1.1}\label{sec:core}
In this section, we identify the pole core of a transversally
Cantor-like $\omega$-limit set and establish restrictions on
the residues of Fuchsian poles using the local developing map. We prove our first main result.

Let $p$ be a Fuchsian pole with residue $\rho$. If $\rho$ is not a negative integer strictly smaller than $-1$, there is a coordinate $z$ centered at $p$ in which
\begin{equation}\label{eq:normal}
 \eta=\frac{\rho}{z}\,dz.
\end{equation}
see \cite{AR}. Put $a:=1+\rho$. The following formulas are used for $a\neq0$; when $a=0$, a developing map on the universal cover is $J(\zeta)=\zeta$.
Using the notation introduced above, we can express the developing map as follows: $J(z)=\frac{z^a}{a}$.
On the universal cover of the punctured coordinate disc, with $z=e^\zeta$, a developing map is
\begin{equation}\label{eq:develop}
 J(\zeta)=\frac{e^{a\zeta}}{a}.
\end{equation}

It is clear that the points $\zeta$ and $\zeta + 2 \pi i$ correspond to exactly the same point under the universal covering map. However, on the universal cover, they are distinct points. Thus, the universal covering allows us to keep track of how many times we wind around $z$. The map $\zeta\mapsto\zeta+2\pi i$ on the universal cover is called a \emph{deck transformation}.
The deck transformation $\zeta\mapsto\zeta+2\pi i$ acts in the developed plane by
\begin{equation}\label{eq:q}
 J(\zeta+2\pi i)=\ qJ(\zeta),
 \qquad q:=e^{2\pi ia}.
\end{equation}

A lifted geodesic is mapped to a Euclidean straight line. In the range $\operatorname{Re}a>0$ considered below, a \emph{critical geodesic} corresponds to a line through the origin; a \emph{noncritical local passage} is a maximal connected component of a geodesic in the punctured coordinate disc whose developed line avoids the origin.

\begin{lemma}\label{lem:complex-self}
Assume
$$
 0<\operatorname{Re}a<\frac12.
$$
For every sufficiently small coordinate disc $\Delta_r=\{|z|<r\}$ there exists $\delta_r>0$ with the following property: every noncritical local geodesic passage in $\Delta_r^*$ whose minimum coordinate modulus is smaller than $\delta_r$ self-intersects before leaving $\Delta_r^*$.
\end{lemma}

\begin{proof}
Let $a=u+iv$. Then $u=\operatorname{Re}a$ and $v=\operatorname{Im}a$. By the hypothesis of the lemma, $u\in\left(0,\frac12\right)$. From $q=e^{2 \pi ia}$ we obtain the following: $q=e^{2 \pi i \cdot(u+iv)}=e^{-2 \pi v} \cdot e^{2 \pi ui}$ and $ q=Re^{i\phi}$.
We get that $R=e^{-2\pi\operatorname{Im}a}>0$ and $\phi=2\pi u\in(0,\pi)$.
Set
\begin{equation}\label{eq:s0}
 s_0:=\frac{R^{-1}-\cos\phi}{\sin\phi},
 \qquad
 w_0:=s_0+i.
\end{equation}
By construction $w_1:=qw_0$ and the point also has imaginary part equal to $1$.

$$
\operatorname{Im}w_1=R(s_0\sin\phi+\cos\phi)=1.
$$
Hence $w_1=s_1+i$, where
$$
s_1=\frac{\cos\phi-R}{\sin\phi}.
$$
Moreover,
$$
s_0-s_1=\frac{R+R^{-1}-2\cos\phi}{\sin\phi}>0.
$$
 Since both points lie in the upper half-plane and multiplication by $q$ rotates by the angle $\phi\in(0,\pi)$, the continuous logarithm on the upper half-plane satisfies
\begin{equation}\label{eq:log-difference}
 \log w_1-\log w_0=\log R+i\phi=2\pi ia.
\end{equation}

Every noncritical developed line can, after multiplication by a nonzero complex constant and an affine reparametrization along the line, be normalized to $\{t+i:t\in\R\}$. Correspondingly, a branch of its inverse image has the form
\begin{equation}\label{eq:zC}
 z_C(t)=C\exp\!\left(\frac1a\log(t+i)\right),
 \qquad C\in\C^*.
\end{equation}
Using \eqref{eq:log-difference}, we obtain $z_C(s_1)
 =z_C(s_0)\exp\!\left(\frac{2\pi ia}{a}\right)
 =z_C(s_0)$.
Thus the full local branch has a self-intersection at the two distinct parameters $s_0$ and $s_1$. It remains to place this self-intersection inside the prescribed coordinate disc. Put:
$$
 F_a(t):=\left|\exp\!\left(\frac1a\log(t+i)\right)\right|.
$$
A direct differentiation gives $\frac{d}{dt}\log F_a(t)
 =\frac{ut-v}{(u^2+v^2)(1+t^2)}$.
Hence $F_a$ has a unique minimum, at $t=v/u$, and tends to $+\infty$ as $|t|\to\infty$. In particular, every sublevel set of $F_a$ is an interval and
$$
 m_a:=\min_{t\in\R}F_a(t)>0.
$$
Let $I_a$ be the compact interval with endpoints $s_0,s_1$, and set
$$
 M_a:=\max_{t\in I_a}F_a(t)<+\infty.
$$
The minimum modulus of \eqref{eq:zC} is $|C|m_a$, while every point between the two self-intersection parameters has modulus at most $|C|M_a$. Therefore, if
$$
 \min_t|z_C(t)|<\delta_r:=r\frac{m_a}{M_a},
$$
the whole subarc between $s_0$ and $s_1$ lies in $\Delta_r^*$ and self-intersects there. Since the relevant sublevel set of $F_a$ is an interval, this subarc belongs to the same connected local passage through $\Delta_r^*$.
\end{proof}

\begin{remark}\label{rem:multiplicity}
The same argument with $q^n$ shows that every positive integer $n$ satisfying
$$
 0<n\operatorname{Re}a<\frac12
$$
produces a forced local self-intersection class. Thus a sufficiently deep passage has at least
$$
 \left\lceil\frac{1}{2\operatorname{Re}a}\right\rceil-1
$$
such classes when $0<\operatorname{Re}a<1/2$.

 Conversely, if $z_C(t_1)=z_C(t_2)$ with $t_1\neq t_2$, then
$$
\log(t_1+i)-\log(t_2+i)=2\pi i k a
$$
for some $k\in\mathbb Z\setminus\{0\}$. Since both arguments
belong to $(0,\pi)$, this implies
$$
0<|k|\operatorname{Re}a<\frac12.
$$
Thus, for $\operatorname{Re}a>0$, the preceding construction
gives exactly the possible deck-intersection classes. In
particular, every noncritical inverse-line branch is injective
when $\operatorname{Re}a\geq1/2$, so the threshold is locally sharp.
\end{remark}

\begin{proof}[Proof of \Cref{thm:core-intro}]

A nonempty compact invariant set contains a minimal compact invariant subset. By \Cref{thm:minimal-AB}, such a subset is a pole, the support of a closed geodesic, or all of $S$. The last alternative is impossible because $W$ has empty interior. The closed-geodesic alternative is also impossible, because \Cref{thm:minimal-AB} would then give that $W$ itself is the support of that closed geodesic, whereas a short transversal to a closed geodesic meets it in a finite set. Thus every minimal subset is a pole, proving in particular that $P_W$ is nonempty.

Conversely, every $p\in P_W$ is a singular invariant point of the associated line field, so $\{p\}$ is a minimal invariant set. The pole set of a meromorphic connection on a compact Riemann surface is finite. Finally, since $\Int W=\varnothing$, every regular geodesic contained in $W$ is distinguished \cite[Proposition~4.1(iii)]{AB}. The support of a closed regular geodesic contained in $W$ would therefore be a minimal invariant subset, which has already been excluded.

 As established in \cite{R}, if a Fuchsian pole $p\in W$ satisfies $\operatorname{Re}\Res_p\nabla\leq-1$, then $W=\{p\}$, impossible here. Suppose then that
$$
 -1<\operatorname{Re}\Res_p\nabla< -\frac12.
$$
In this range the pole is nonresonant, so the normal form \eqref{eq:normal} is available, and $a=1+\Res_p\nabla$ satisfies
$$
 0<\operatorname{Re}a<\frac12.
$$

Because $p\in W=\overline{\supp\sigma}$ and $W$ is nontrivial, the simple geodesic has noncritical passages approaching $p$ arbitrarily closely. Indeed, a passage critical in the forward direction would terminate the maximal geodesic at $p$ and force $W=\{p\}$; a critical backward passage can occur only as an initial end and is irrelevant after deleting an initial segment. Applying \Cref{lem:complex-self} to a sufficiently deep noncritical passage produces a self-intersection, contradicting simplicity.
\end{proof}

 The residue argument only uses eventual simplicity and $W\neq\{p\}$. Hence the same bound holds at every Fuchsian pole in any nontrivial $\omega$-limit set of an eventually simple geodesic.

\begin{corollary}\label{cor:poles-not-isolated}
No pole in $P_W$ is isolated in $W$.
\end{corollary}

It is not difficult to see that a nondegenerate connected compact set has no isolated point.

\begin{corollary}\label{cor:necessary-general}
For a Fuchsian meromorphic connection admitting a transversally Cantor-like $\omega$-limit set $W$, the finite pole core $P_W$ is nonempty and satisfies
$$
 \operatorname{Re}\Res_p\nabla\geq-\frac12
 \qquad\text{for every }p\in P_W.
$$
\end{corollary}
\section{Proof of Theorem 1.2}\label{sec:real period}

 We prove \Cref{thm:flat-intro} using regular return bands and transversal gaps.

\subsection{\texorpdfstring{ Real periods and regular crossing states}{Real periods and regular crossing states}}\label{sec:gates}

 Let $\nabla$ be Fuchsian with real periods, and let $g$ be adapted to $\nabla$.
In flat coordinates,
\begin{equation}\label{eq:unitary-transition}
 w_\beta=e^{i\theta_{\alpha\beta}}w_\alpha+c_{\alpha\beta}.
\end{equation}
Assume that a simple tail $\sigma$ has a transversally Cantor-like limit set $W$.
Normalize its parameter so that
$$
 |\sigma'|_g\equiv1,\qquad
 W=\omega(\sigma)=\overline{\supp\sigma},\qquad \Int W=\varnothing.
$$
By \Cref{thm:core-intro}, $W$ contains no closed regular geodesic.

Choose a short compact geodesic transversal $I$ through $\supp\sigma$, perpendicular
to the distinguished line at its central point, with
$$
 \partial I\cap W=\varnothing,\qquad C:=I\cap W\ \text{a Cantor set}.
$$
Shrink $I$ so that it has a pole-free Euclidean collar and the distinguished lines
make an angle $|\theta(s)|<\varepsilon_0<\pi/8$ with its normal.

 Set
$$
 \mathcal I:=I_+\bigsqcup I_-,\qquad \pi:\mathcal I\longrightarrow I.
$$
On $I_\pm$, let $u_\pm$ be the forward unit normal and choose the unit tangent
$e_\pm$ so that $(u_\pm,e_\pm)$ is positively oriented. Use signed flat arclength
in the direction $e_\pm$. The distinguished forward unit tangent is
\begin{equation}\label{eq:gate-angle}
 v(s)=\cos\theta(s)\,u_\pm+\sin\theta(s)\,e_\pm,
 \qquad |\theta(s)|<\varepsilon_0,
 \qquad \theta\in C(\pi^{-1}(C),\R).
\end{equation}
Record consecutive oriented crossings as
$$
 x_n=(\sigma(t_n),\varepsilon_n)\in\mathcal I,
 \qquad \varepsilon_n\in\{+,-\},\qquad t_0<t_1<\cdots,
$$
and put
\begin{equation}\label{eq:phase-K}
 K:=\bigcap_{N\geq0}\overline{\{x_n:n\geq N\}}\subset\mathcal I.
\end{equation}

\begin{lemma}\label{lem:directions}\label{lem:crossing-states}
The set $K$ is compact and
$$
 t_n\longrightarrow+\infty,\qquad
 \pi(K)=C,\qquad \#K=\infty,\qquad \Int_{\mathcal I}K=\varnothing.
$$
Moreover, $\theta|_K\in C(K,\R)$.
\end{lemma}

\begin{proof}
Uniform transversality turns arbitrarily late visits to $q\in C$ into nearby
crossings. Selecting a crossing direction gives $\pi(K)=C$.
Since $C\Subset\Int I$, a unit-speed return must leave a fixed smaller collar;
hence, for some $\delta>0$,
$$
 t_{n+1}-t_n\geq\delta,\qquad t_n\geq t_0+n\delta\longrightarrow+\infty.
$$
By \eqref{eq:phase-K}, $K$ is compact. Its projection is the Cantor set $C$, so
$\#K=\infty$ and $\Int_{\mathcal I}K=\varnothing$.
Continuity of $\theta|_K$ follows from \eqref{eq:gate-angle}.
\end{proof}

 We use $K$ without assuming $x_n\in K$, $K\cap I_\pm$ perfect, or finitely many exact crossing directions.

\subsection{Finite Euclidean return bands}\label{sec:bands}

 Blow up the poles and slit along $I$, marking both slit tips. The resulting
compact surface $M$ contains pairwise disjoint properly embedded return arcs
$$
 a_n:=\sigma([t_n,t_{n+1}]),\qquad a_n\cap a_m=\varnothing\quad(n\neq m).
$$
At a shared crossing, incoming and outgoing endpoints lie on opposite physical
slit sides.

Write $a\sim b$ if their oriented source and target copies and physical slit
sides agree, and they are relatively isotopic in $M$. Endpoints may slide within
the prescribed open sides, without passing a slit tip. If both endpoints lie
on one physical side, retain their source--target order. Put
$$
 \mathcal T:=\{[a_n]_\sim:n\geq0\}.
$$

\begin{lemma}\label{lem:finite-types}
Only finitely many topological types occur among the return arcs $a_n$.
\end{lemma}

\begin{proof}
A finite ideal triangulation bounds the number $N(M)$ of pairwise
disjoint, nonisotopic essential arcs relative to the boundary marks.
An inessential return arc would bound a pole-free geodesic bigon. If its angles
are $\alpha,\beta>0$, flat Gauss--Bonnet gives
$$
 (\pi-\alpha)+(\pi-\beta)=2\pi
 \quad\Longrightarrow\quad \alpha+\beta=0,
$$
a contradiction. Boundary-parallel arcs enclosing marks are essential relative
to those marks. Since there are only $L<\infty$ endpoint flags,
$$
 \#\mathcal T\leq LN(M)<\infty.
$$
\end{proof}

After deleting finitely many initial returns, we may assume that every remaining type occurs infinitely often.

\begin{lemma}\label{lem:rectangle}\label{lem:quadrilateral}
Two distinct arcs $a\sim b$ and their intervening gate segments bound a pole-free
product disk $Q$, isometric to a convex Euclidean quadrilateral $P$.
Distinguished geodesics entering between the source endpoints exit between the
target endpoints, and conversely. Their unit-speed transit times satisfy
$$
 0<T_Q(x)\leq\operatorname{diam}P.
$$
\end{lemma}

\begin{proof}
Relative isotopy gives a product disk with
$$
 \partial Q=a\cup c_R\cup b\cup c_D,\qquad Q\cap\Sing(\nabla)=\varnothing,
$$
containing no slit tip or other boundary component. Its developing map is a
local Euclidean isometry. All four sides develop to straight segments; their
interior angles $\alpha_j$ satisfy
$$
 0<\alpha_j<\pi,\qquad \sum_{j=1}^4(\pi-\alpha_j)=2\pi.
$$
Hence the developed boundary is a simple convex quadrilateral $\partial P$.
The boundary degree is one, so $\operatorname{dev}:Q\to P$ is an isometry.

A distinguished plaque cannot cross $a$ or $b$, nor meet the straight source
side twice. It therefore exits through $c_R$, with transit time at most
$\operatorname{diam}P$. Reverse time for the converse.
\end{proof}

 For $\tau\in\mathcal T$, define, on the prescribed oriented copies,
$$
 \begin{aligned}
 S_\tau&:=\{x_n:[a_n]_\sim=\tau\},&
 E_\tau&:=\{x_{n+1}:[a_n]_\sim=\tau\},\\
 D_\tau&:=\Int\operatorname{conv}S_\tau,&
 R_\tau&:=\Int\operatorname{conv}E_\tau.
 \end{aligned}
$$
These are nondegenerate open hulls. By \Cref{lem:rectangle}, a crossing between
two source points of type $\tau$ has the same type; the analogous statement
holds backward.

\begin{proposition}\label{prop:partial-IET}\label{prop:ordered-returns}
For $\tau\neq\kappa$,
$$
 D_\tau\cap D_\kappa=R_\tau\cap R_\kappa=\varnothing.
$$
Distinguished returns induce increasing homeomorphisms
\begin{equation}\label{eq:K-return}
 F_\tau:K\cap D_\tau\longrightarrow K\cap R_\tau.
\end{equation}
Put
$$
 D:=\bigcup_\tau D_\tau,\quad R:=\bigcup_\tau R_\tau,\quad
 E_D:=\bigcup_\tau\partial D_\tau,\quad E_R:=\bigcup_\tau\partial R_\tau.
$$
Then
\begin{equation}\label{eq:K-hull-cover}
 K\subset\overline D\cap\overline R,\qquad
 K\setminus D\subset E_D,\qquad K\setminus R\subset E_R.
\end{equation}
In particular, both exceptional sets are finite.
\end{proposition}

\begin{proof}
Overlapping source hulls would place a crossing of one type between two of
another, contrary to trapping in \Cref{lem:rectangle}. Reverse time for the
range hulls.

Given $x\in K\cap D_\tau$, enclose it between two observed source points of type
$\tau$, bounding a fixed quadrilateral $Q$. For late crossings $x_{n_j}\to x$,
trapping and continuous dependence give
$$
 0<t_{n_j+1}-t_{n_j}\leq\operatorname{diam}P,\qquad
 F_\tau(x)=\lim_jx_{n_j+1}\in K\cap R_\tau.
$$
Conversely, backward trapping at $y\in K\cap R_\tau$ gives
$$
 F_\tau^{-1}(y)=\lim_jx_{n_j-1}\in K\cap D_\tau.
$$
The bounded transit times ensure that these predecessors are late crossings.
Continuous dependence gives continuity of $F_\tau^{\pm1}$; noncrossing and the
chosen frames give increasing order.

Except for at most two extrema per type, every observed crossing belongs to
its open hull. Simplicity allows at most one visit to each extremum.
Taking tail closures proves \eqref{eq:K-hull-cover}.
\end{proof}

\subsection{\texorpdfstring{Gap chains and exclusion of Cantor recurrence}{Gap chains and exclusion of Cantor recurrence}}\label{sec:IET}

 Let $\mathcal G$ be the bounded complementary intervals of $K$ in the two
oriented copies. Since $K$ is infinite, compact and nowhere dense,
$$
 \#\mathcal G=\infty,\qquad
 G=(a,b)\in\mathcal G\ \Longrightarrow\ a,b\in K,\quad |G|:=b-a>0.
$$

\begin{lemma}\label{lem:completion}\label{lem:gap-permutation}
The ordered returns define a partial permutation
$\Phi:\operatorname{Dom}\Phi\to\operatorname{Ran}\Phi$ of $\mathcal G$, with
$$
 \#(\mathcal G\setminus\operatorname{Dom}\Phi)
 +\#(\mathcal G\setminus\operatorname{Ran}\Phi)<\infty.
$$
\end{lemma}

\begin{proof}
For $G=(a,b)$ with $a,b\in D_\tau$, set
$$
 \Phi(G):=(F_\tau(a),F_\tau(b)).
$$
By \eqref{eq:K-return}, an intermediate target point of $K$ would have an
intermediate source point in $K$. Thus $\Phi(G)$ is a gap. Disjoint hulls give
uniqueness; backward returns give the inverse.

By \eqref{eq:K-hull-cover}, a missing forward arrow forces
$\overline G\cap E_D\neq\varnothing$. A point belongs to at most two gap
closures, so
$$
 \#(\mathcal G\setminus\operatorname{Dom}\Phi)\leq2\#E_D,\qquad
 \#(\mathcal G\setminus\operatorname{Ran}\Phi)\leq2\#E_R.
$$
\end{proof}

\begin{lemma}\label{lem:IET-dichotomy}\label{lem:angle-width}
For $a<b$ in $K\cap D_\tau$,
\begin{equation}\label{eq:angle-difference-invariant}
 \theta(F_\tau(b))-\theta(F_\tau(a))=\theta(b)-\theta(a).
\end{equation}
If $\theta(a)=\theta(b)$, the developed boundary geodesics are parallel and
\begin{equation}\label{eq:normal-gap-width}
 h(a,b):=(b-a)\cos\theta(a)>0,\qquad
 h(F_\tau(a),F_\tau(b))=h(a,b).
\end{equation}
\end{lemma}

\begin{proof}
Closing return paths by gate segments fixes their homotopy class and hence their
unitary holonomy. Including the fixed source and target frames, there is a
rotation $\beta_\tau$ such that
$$
 \theta(F_\tau(s))\equiv\theta(s)+\beta_\tau\pmod{2\pi}
 \qquad(s\in K\cap D_\tau).
$$
Since $|\theta(s)|<\varepsilon_0<\pi/8$, subtraction has no $2\pi$ ambiguity, proving
\eqref{eq:angle-difference-invariant}.

If $\theta(a)=\theta(b)$, develop a quadrilateral enclosing $a,b$.
The parallel boundary lines have perpendicular distance
$$
 (b-a)\cos\theta(a)
 =(F_\tau(b)-F_\tau(a))\cos\theta(F_\tau(a)).
$$
Unitary transitions preserve this distance, proving \eqref{eq:normal-gap-width}.
\end{proof}

\begin{proof}[Proof of \Cref{thm:flat-intro}]
A gap cycle $\Phi^m(G)=G$ returns both ordered endpoint states to themselves.
For an endpoint geodesic $\gamma$, this gives, for some $L>0$,
$$
 (\gamma(L),\gamma'(L))=(\gamma(0),\gamma'(0)),\qquad
 \supp\gamma\subset W.
$$
Uniqueness yields a closed regular geodesic, contradicting
\Cref{thm:core-intro}. Thus $\Phi$ has no cycles.

By \Cref{lem:gap-permutation}, only finitely many arrows are missing.
Since $\#\mathcal G=\infty$, there is a forward or backward infinite chain of
distinct gaps $G_n=(a_n,b_n)$. Put $\ell_n:=b_n-a_n$. Disjointness gives
$$
 \sum_{n\geq0}\ell_n\leq2|I|<\infty
 \quad\Longrightarrow\quad\ell_n\longrightarrow0.
$$
By \eqref{eq:angle-difference-invariant} and uniform continuity of $\theta|_K$,
$$
 \Delta_n:=\theta(b_n)-\theta(a_n)=\Delta_0,\qquad
 \Delta_n\longrightarrow0.
$$
Hence $\Delta_n=0$ for every $n$. Applying \eqref{eq:normal-gap-width} in either
direction now gives
$$
 0<h(a_0,b_0)=h(a_n,b_n)
   =\ell_n\cos\theta(a_n)\leq\ell_n\longrightarrow0,
$$
a contradiction.
\end{proof}

\begin{corollary}\label{cor:nonunitary}
Let $\nabla$ be Fuchsian and suppose it admits a transversally Cantor-like $\omega$-limit set. Then its modulus character
$$
 \ell_\nabla:\pi_1(S^0)\longrightarrow\R,
 \qquad
 \ell_\nabla(\gamma):=\log|\Hol_\nabla(\gamma)|,
$$
is nonzero.
\end{corollary}

\begin{proof}
If $\ell_\nabla\equiv0$, then
$$
 |\Hol_\nabla(\gamma)|=e^{\ell_\nabla(\gamma)}=1
 \quad(\gamma\in\pi_1(S^0))
 \quad\Longrightarrow\quad
 \Hol_\nabla(\pi_1(S^0))\subset\Sone.
$$
Thus $\nabla$ has real periods, contrary to \Cref{thm:flat-intro}.
\end{proof}

\section{A strengthened Poincar\'e--Bendixson theorem}\label{sec:PB}

The general Poincar\'e--Bendixson theorem of \cite{AB, AT}, with the missing Cantor-like case included as in \cite{AR}, lists poles, closed geodesics, closed-geodesic supports, boundary graphs of saddle connections, Cantor-like sets, the whole surface, sets with nonempty interior and saddle-connection boundary, and infinite self-intersection. For real periods, closed geodesics are periodic; moreover, \cite{AR} excludes a periodic-geodesic support and a boundary graph of saddle connections as the $\omega$-limit set of a simple nonperiodic geodesic. Together with \Cref{thm:flat-intro}, this gives the following.

\begin{corollary}\label{cor:PB}
Let $\nabla$ be a Fuchsian meromorphic connection with real periods on a compact Riemann surface $S$. Let
$$
 \sigma:[0,\varepsilon)\longrightarrow S^0
$$
be a maximal geodesic with at most finitely many transverse self-intersections. After deleting an initial segment, exactly one of the following occurs:
\begin{enumerate}
\item $\omega(\sigma)$ is a pole of $\nabla$;
\item $\sigma$ is periodic;
\item $\omega(\sigma)$ has nonempty interior and nonempty boundary, and every boundary component is a graph of saddle connections without spikes and containing at least one pole;
\item $\omega(\sigma)=S$.
\end{enumerate}
\end{corollary}

\begin{proof}
 If $\sigma$ is periodic, we are in case~(ii). Otherwise, a tangential return would imply periodicity, by uniqueness for the geodesic equation and constancy of the nonzero adapted speed. Delete the portion preceding the last self-intersection; the remaining geodesic is simple and has the same $\omega$-limit set. Apply the general classification. For a nonperiodic geodesic, the closed-geodesic-support, boundary-graph, and Cantor-like alternatives are excluded by the results just recalled and \Cref{thm:flat-intro}. The four displayed cases are the remaining possibilities.
\end{proof}

\subsection{Affine returns and the Disco surface}

A finite affine interval exchange is a piecewise affine bijection
with branches
\begin{equation}\label{eq:AIET}
 T(x)=\lambda_jx+c_j,\qquad\lambda_j>0.
\end{equation}
When all slopes are one, this is an ordinary interval exchange.
For general slopes, wandering intervals can occur, and their
complement can support exceptional Cantor recurrence.
This phenomenon was constructed in \cite{CG} and developed further
in the self-similar and typical settings of \cite{BHM,MMY}.

An affine suspension glues transversal intervals of rectangles by
the maps in \eqref{eq:AIET}, using complex affine gluings that preserve
the vertical direction. Away from singular vertices it has coordinates
with transition maps
\begin{equation}\label{eq:dilation-transition}
 z_\beta=a_{\alpha\beta}z_\alpha+b_{\alpha\beta},
 \qquad a_{\alpha\beta}>0.
\end{equation}
Such a dilation surface carries a directional foliation whose
first-return map is $T$; see \cite{BFG, DFG}.
The local forms $\eta_\alpha=0$ define a holomorphic connection on
the regular part. Standard Fuchsian affine singularities allow its
meromorphic extension to the compact Riemann surface \cite{NST,R}.

\begin{example}\label{ex:disco}
The Disco surface of \cite{BFG} is a compact dilation surface of genus
two with two Euclidean cone points of angle $4\pi$.
Near each point, a coordinate $z$ has a developing map $J=z^2/2$,
up to affine normalization. Its induced connection has
$$
 \eta=d\log J'=\frac{dz}{z},
$$
so the two Fuchsian residues are
\begin{equation}\label{eq:disco-residues}
 \rho_1=\rho_2=1.
\end{equation}
With the oriented symplectic basis used in { \cite[Section~2.6 of the earlier preprint]{BFG}},
the linear developing holonomy satisfies
\begin{equation}\label{eq:disco-holonomy}
 \Hol(a_1)=\Hol(b_2)=2,\qquad
 \Hol(a_2)=\Hol(b_1)=\frac12.
\end{equation}
In particular, this connection has real residues and fails the
real-period condition because its linear holonomy is nonunitary.

 For a Cantor direction supplied by \cite[Theorem~2]{BFG}, let $T$ be the
return map and $C$ its common Cantor $\omega$-limit set. If $\mathcal D$ denotes
the finite set of discontinuities of $T$ and $T^{-1}$, choose
$$
 x\in C\setminus\bigcup_{n\in\mathbb Z}T^n(\mathcal D),\qquad
 \omega_T(x)=C,\qquad x\in\omega_T(x).
$$
Its orbit is nonperiodic and avoids singular return paths. The suspended
geodesic is simple, is contained in its own $\omega$-limit set, and has Cantor
transverse sections. Hence it satisfies \Cref{def:cantor}.
\end{example}


\begin{thebibliography}{99}

\bibitem{AB}
M.~Abate and F.~Bianchi,
\emph{A Poincar\'e--Bendixson theorem for meromorphic connections on Riemann surfaces},
Math. Z. \textbf{282} (2016), 247--272.

\bibitem{AR}
M.~Abate and K.~Rakhimov,
\emph{Dynamics of Fuchsian meromorphic connections with real periods},
 Ann. Global Anal. Geom. \textbf{69} (2026), Article~9.
\href{https://doi.org/10.1007/s10455-026-10033-8}{doi:10.1007/s10455-026-10033-8}.

\bibitem{AT}
M.~Abate and F.~Tovena,
\emph{Poincar\'e--Bendixson theorems for meromorphic connections and homogeneous vector fields},
J. Differential Equations \textbf{251} (2011), 2612--2684.

\bibitem{BFG}
A.~Boulanger, C.~Fougeron and S.~Ghazouani,
\emph{Cascades in the dynamics of affine interval exchange transformations},
Ergodic Theory Dynam. Systems \textbf{40} (2020), no.~8, 2073--2097.

\bibitem{BHM}
X.~Bressaud, P.~Hubert and A.~Maass,
\emph{Persistence of wandering intervals in self-similar affine interval exchange transformations},
Ergodic Theory Dynam. Systems \textbf{30} (2010), no.~3, 665--686.

\bibitem{CG}
R.~Camelier and C.~Gutierrez,
\emph{Affine interval exchange transformations with wandering intervals},
Ergodic Theory Dynam. Systems \textbf{17} (1997), no.~6, 1315--1338.

\bibitem{Dah}
F.~Dahmani,
\emph{Groups of interval exchange transformations},
Winter Braids Lecture Notes \textbf{6} (2019), Course no.~I, 1--22.

\bibitem{DFG}
E.~Duryev, C.~Fougeron and S.~Ghazouani,
\emph{Dilation surfaces and their Veech groups},
J. Mod. Dyn. \textbf{14} (2019), 121--151.

\bibitem{Hou}
J.~Hounie,
\emph{Minimal sets of families of vector fields on compact surfaces},
J. Differential Geom. \textbf{16} (1981), 739--744.

\bibitem{Ima}
H.~Imanishi,
\emph{On codimension one foliations defined by closed one-forms with singularities},
J. Math. Kyoto Univ. \textbf{19} (1979), no.~2, 285--291.

\bibitem{MMY}
S.~Marmi, P.~Moussa and J.-C.~Yoccoz,
\emph{Affine interval exchange maps with a wandering interval},
Proc. Lond. Math. Soc. (3) \textbf{100} (2010), no.~3, 639--669.

\bibitem{NST}
D.~Novikov, B.~Shapiro and G.~Tahar, 
\emph{On limit sets for geodesics of meromorphic connections},
J. Dyn. Control Syst. \textbf{29} (2023), 55--70.
\href{https://doi.org/10.1007/s10883-021-09584-9}{doi:10.1007/s10883-021-09584-9}.

\bibitem{R}
K.~Rakhimov,
\emph{Flat structure of meromorphic connections on Riemann surfaces},
J. Geom. Anal. \textbf{35} (2025), Article~368.


\end{thebibliography}
\end{document}